\documentclass[leqno,11pt,a4paper]{article}
\usepackage{pdfsync}
\usepackage{amsmath}
\usepackage{amsfonts}
\usepackage{amssymb}
\usepackage{amsthm}
\usepackage{inputenc}
\usepackage{enumitem,fullpage}
\usepackage{hyperref}
\usepackage{graphicx,multirow}
\usepackage{xcolor}
\usepackage{mwe}
\usepackage{tikz}
\usepackage{tikz-cd}
\usepackage[margin=1in, includeheadfoot]{geometry} 
\hypersetup{
 colorlinks,
 linkcolor={blue!50!black},
 citecolor={blue!50!black},
 urlcolor={blue!80!black}
}
\usetikzlibrary{circuits.logic.US,circuits.logic.IEC,fit}

\usepackage{arydshln}
\usepackage{multirow}
\usepackage[all]{xy}
\usepackage{mleftright}
\usepackage{url}

\usepackage{pdflscape}
\usepackage{longtable}
\usepackage{booktabs}
\usepackage{colortbl}
\newcommand{\evnrow}{\rowcolor[gray]{0.95}}
\newcommand{\oddrow}{}

\usepackage{array,etoolbox}
\preto\tabular{\setcounter{magicrownumbers}{0}}
\newcounter{magicrownumbers}
\newcommand\rownumber{\stepcounter{magicrownumbers}\arabic{magicrownumbers}}
\setlist[enumerate]{labelsep=*, leftmargin=1.5pc}
\setlist[enumerate]{label=\normalfont(\roman*), ref=\roman*}

\newcommand{\C}{\mathbb C}
\renewcommand{\P}{\mathbb P}
\newcommand{\Af}{\mathbb A}
\newcommand{\Q}{\mathbb Q}
\newcommand{\M}{\mathbb M}

\newcommand{\Z}{\mathbb Z}

\newcommand{\A}{\mathbb A}
\newcommand{\V}{\mathbb V}
\newcommand{\cB}{\mathcal{B}}
\newcommand{\cO}{\mathcal{O}}

\newcommand{\cC}{\mathcal{C}}
\newcommand{\I}{\mathbb{I}}

\def\wP{\mathrm{w}\mathcal{P}}

\newcommand{\orb}{{\mathrm{orb}}}

\newcommand{\into}{\hookrightarrow}
\def\sm{\mathrm{smooth}}

\DeclareMathOperator{\Gr}{Gr}

\DeclareMathOperator{\GL}{{GL}}
\DeclareMathOperator{\Hom}{Hom}

\newcommand{\PxP}{\P^2\times\P^2}
\newcommand{\PI}{\P^1\times\P^1 \times\P^1}

\newcommand{\Oh}{\mathcal{O}}

\newtheorem{thm}{Theorem}

\newtheorem{prop}[thm]{Proposition}

\theoremstyle{definition}
\newtheorem{dfn}[thm]{Definition}
\theoremstyle{remark}
\newtheorem{rmk}[thm]{Remark}
\theoremstyle{remark}
\newtheorem{eg}[thm]{Example}

\numberwithin{equation}{section}
\numberwithin{thm}{section}

\newcommand{\QED}{\ifhmode\unskip\nobreak\fi\quad\ensuremath{\mathrm{QED}}}

\title{Orbifold del Pezzo surfaces  in $\P^1\times\P^1 \times \P^1$ format}
\author{ Muhammad Imran Qureshi}
\date{}

\begin{document}

\maketitle

\begin{abstract}
    We introduce  the notion of biregular index models    (infinite series)  of  orbifold del Pezzo surfaces having their    (sub) anti-canonical  embeddings  in some weighted projective space   \(\P(w_i) \). We construct such models of
del Pezzo surface with  embeddings in \(\P^6(w_i)\) containing at worst
rigid orbifold points. The equations describing their images under    their
(sub) anti-canonical embeddings can be computed by using   equations of the Segre embedding of \(\PI\) in
\(\P^7\); giving  rise to codimension 4 Gorenstein varieties.  We also compute a     formula for the Hilbert series of a general  weighted \(\PI\into\P^7(w_i)\) variety
that plays a pivotal role in  these constructions.

\end{abstract}

%
%%%%%%%%%%%%%%%%%%%%%%%%%%%%%%%%%%%%
\section{Introduction}
\subsection{Preliminaries}
An algebraic surface $X$ is a \emph{del Pezzo surface} if the anti-canonical class $-K_X$ is ample. The \emph{Fano index} \(I\) of \(X\) is the largest positive integer $I$ such that $-K_X = I\cdot D$ for some divisor $D $ in the class group of $X$. One can study   various properties of a del Pezzo surface \(X\) by examining its  images under  (sub) anti-canonical embeddings
in some  weighted projective space \(\P^k(w_i)\). If the image can be described by using a simple set of equations,
like minors or Pfaffians of some matrix then we say that the surface \(X\)
is given by an equation format.   If the affine cone  \(\widetilde
X \subset \A^{k+1}\) of the embedding  is smooth outside the origin then we call \(X\) to be quasismooth. If the codimension of \(X\into \P^k(w_i)\)
is \(c\) then we call  \(X\) to be  wellformed if   does not contain
  any codimension \(c+1\)
orbifold locus of \(\P^k(w_i)\)  .  
 
An orbifold point    of type  $\frac{1}{r}(a,b)$ on a surface  is the quotient \(\pi:\Af^ 2\to \Af^2/\mu_r\) given by \[\mu_r \ni\epsilon:(x_1,x_2)\mapsto (\epsilon^{a}x_1,\epsilon^{b}x_2).  \]  It is called an isolated orbifold point if   $r$ is relatively prime to both  $a$ and $b$. A  del Pezzo surface containing at worst orbifold points  is called a log or an orbifold del Pezzo surface (odP). 

A rigid del Pezzo surface is a del Pezzo surface whose singularities are not smoothable by  $\Q$-Gorenstein deformations ~\cite{Kollar-SB}. We use the following characterization of rigid orbifold points to check for the rigidity of the given  del Pezzo surface.  

\begin{dfn}~\cite{AK-SC}
\label{dfn:T_R}
Let  $Q = \frac{1}{r}(a,b)$ be an orbifold point and   $d = \text{gcd}(a+b,r)$, $k = (a+b)/d$ and $m=r/d$. Then we can present \(Q\) as    $\frac{1}{dm}(1,mk-1)$ and if $d<m$ then it is called a (rigid) \emph{R}-singularity.
\end{dfn}  
We recall the notion of  biregular models of orbifold del Pezzo
surfaces,  introduced in~\cite{QMOC}
\begin{dfn}
A \emph{biregular model of orbifold del Pezzo surfaces} \(\M\) is  an infinite series of    del Pezzo surfaces satisfying the following conditions.  
\begin{enumerate}

\item  A family of  orbifold  del Pezzo surfaces exists for each  parameter $r(n)$ for all  positive integers $n,   $ where \(r\) is  a  linear function of \(n\). 
\item The (sub) anti-canonical embedding of each family of surfaces in the
weighted projective space \(\P(w_i) \)  has at least one weight equal to   $r$, and  rest of the weights  $w_i$  and the volume  $(-K_X)^2$ are  functions of $r$.  

\end{enumerate}
\end{dfn}
A biregular model \(\M\) is called rigid, if each  del Pezzo surface in
\(\M\) contains at worst rigid orbifold points. We call a biregular model
\(\M\) to be a \textit{weight model} if each odP \(X_r\into \P(w_i)\) in \(\M\) has the  fixed Fano index and the weights \(w_is\)
 are parameterized  by the set
of positive integers.
A series of  hypersurfaces \(X_{2r}\into\P(1,1,r,r) \) of  Fano index \(2\) with two singular points of type \(\frac1r(1,1)\) for \(r=n+4, n\in\Z_+\) is an example of a rigid biregular weight model. In \cite{QMOC}, we proved
the existence of what we  call,  rigid \textit{biregular weight models} in
two  equation formats: a codimension   three  format defined by the maximal Pfaffians of \(5\times 5\)
  skew-symmetric matrix and a codimension four  format defined by \(2\times
  2\) minors of  size 3 matrix. These codimension 3 and 4 formats describe the equations   
   of Grassmannian  \(\Gr(2,5)\) in its Pl\"ucker embedding and of  \(\PxP\)
in its Segre embedding respectively, so we call them \(\Gr(2,5)\) and \(\PxP\)
formats. 
%%%%%%%%%%%%%%%%%%%%%%%%%%%%%%%%%%%%
\subsection{Aim and results}

The  aim of this paper is to   construct  a new type of biregular
models, called\emph{ biregular index models}.     
A \textit{biregular index model} is  an infinite series of  orbifold del Pezzo surfaces where, for each surface,  both  the  Fano index and the set of weights of \(\P^n(w_i),\)   are parameterized by the set of positive integers. 
 A simple example of an index model is the series of surfaces \(\P(1,1,r)\) of index \(r+2\) having an orbifold point of type \(\frac1r(1,1)\) for each \(r=n+1, n \in \Z_+\).

 We construct biregular models of \emph{orbifold rigid del Pezzo} surfaces of codimension 4 in \(\P^6(w_i)\) whose equations can be induced from those of the Segre embedding \(\PI\into \P^7\), i.e. each family of surfaces can be obtained as quasilinear section(s) of some ambient weighted \(\PI\) variety or some projective cone over it ~\cite{QS}. Our main aim is to construct biregular index models, however, we also study and prove the existence of  new weight models in \(\PI\)
format.

\begin{thm} \label{Th:Index}(\textbf{Index models}) There exist at least five  biregular  index  models  of wellformed and quasismooth rigid  del Pezzo surfaces in codimension 4; such that  their images under their sub anti-canonical   embeddings in   \(\P^6(a,\ldots,g)\)  can be described in terms of the equations of the   Segre embedding of \(\PI \) in \(\P^7\).  

 \renewcommand*{\arraystretch}{1.2}
\begin{longtable}{>{\hspace{0.5em}}llccccr<{\hspace{0.5em}}}
\caption{Biregular index models of rigid orbifold del  Pezzo surface in  $\PI $ format of Fano index \(r\). } \label{Index-Models}\\
\toprule
\multicolumn{1}{c}{Model}&\multicolumn{1}{c}{ Parameter  \ }&\multicolumn{1}{c}{Weight Cube}&\multicolumn{1}{c}{\textrm{WPS} $\&\;\;\mathcal{B}$}&\multicolumn{1}{c}{$-K^2$}&\multicolumn{1}{c}{$h^0(-K)$}\\
\cmidrule(lr){1-1}\cmidrule(lr){2-2}\cmidrule(lr){3-3}\cmidrule(lr){4-4}\cmidrule(lr){5-6}
\endfirsthead
\multicolumn{7}{l}{\vspace{-0.25em}\scriptsize\emph{\tablename\ \thetable{} continued from previous page}}\\
\midrule
\endhead
\multicolumn{7}{r}{\scriptsize\emph{Continued on next page}}\\
\endfoot
\bottomrule
\endlastfoot

\evnrow $\M\I_{1} $&$ 
  \begin{matrix}
    \boldsymbol {r=4n+1}\\p=2n,q=r-p\\m=2r-1,s=r+p \\w=(0,p:0,p: 1,r)
  \end{matrix}$&
 \begin{tikzpicture}[
  line join=round,
  y={(-0.86cm,0.36cm)},x={(.6cm,0.36cm)}, z={(0cm,.8cm)},
  arr/.style={line cap=round,shorten <= 1pt},baseline={([yshift=-.5em]current bounding box.center)}
]
\def\Side{1}
\coordinate (A1) at (0,0,0);\coordinate (A2) at (0,\Side,0);
\coordinate (A3) at (\Side,\Side,0);\coordinate (A4) at (\Side,0,0);\coordinate (B1) at (0,0,\Side);\coordinate (B2) at (0,\Side,\Side);\coordinate (B3) at (\Side,\Side,\Side);\coordinate (B4) at (\Side,0,\Side);

\draw[thin] (A2) -- (A1) -- (A4);\draw[thin] (B2) -- (B1) -- (B4) -- (B3) -- cycle;\draw[thin] (A1) -- (B1);\draw[thin] (A2) -- (B2);
\draw[thin] (A4) -- (B4);\draw[thin] (A2) -- (A3);
\draw[thin] (A3) -- (B3);\draw[thin] (A3) -- (A4);

\path[arr] 
  (A1) edge (A2)
  (B2) edge (A2)
  (B1) edge (B2)
  (B1) edge (A1)
  (B4) edge (A4)
  (B3) edge (A3)
  (B4) edge (B3)
  (A4) edge (A3);

\node[below] at (A1) {$1$};
\node[below] at (A2) {$q$};
\node[below] at (A3) {$\boldsymbol r$};
\node[below] at (A4) {$q$};
\node[above] at (B1) {$r$};
\node[above] at (B2) {$s$};
\node[above] at (B3) {$m$};
\node[above] at (B4) {$\boldsymbol s$};

\end{tikzpicture} &
$  \begin{matrix} \P(1,q^2, r,s^2, m)\\[2mm]   
2\times \dfrac1s(1,p-1),\dfrac{1}{m}(1,1)\\2\times \frac1r(p,p+1)
\end{matrix}$&$ \frac{24 r^3}{6 r^3+r^2-4 r+1}       $&$ 3$ \\

\oddrow $\M\I_{2} $&$ 
  \begin{matrix}
    \boldsymbol {r=2n+3},\\p=r-4,q=r-2\\s=r+2,t=r+4 \\w=(0,2:0,2: p,r)
  \end{matrix}$&
 \begin{tikzpicture}[
  line join=round,
  y={(-0.86cm,0.36cm)},x={(.6cm,0.36cm)}, z={(0cm,.8cm)},
  arr/.style={line cap=round,shorten <= 1pt},baseline={([yshift=-.5em]current bounding box.center)}
]
\def\Side{1}
\coordinate (A1) at (0,0,0);\coordinate (A2) at (0,\Side,0);
\coordinate (A3) at (\Side,\Side,0);\coordinate (A4) at (\Side,0,0);\coordinate (B1) at (0,0,\Side);\coordinate (B2) at (0,\Side,\Side);\coordinate (B3) at (\Side,\Side,\Side);\coordinate (B4) at (\Side,0,\Side);

\draw[thin] (A2) -- (A1) -- (A4);\draw[thin] (B2) -- (B1) -- (B4) -- (B3) -- cycle;\draw[thin] (A1) -- (B1);\draw[thin] (A2) -- (B2);
\draw[thin] (A4) -- (B4);\draw[thin] (A2) -- (A3);
\draw[thin] (A3) -- (B3);\draw[thin] (A3) -- (A4);

\path[arr] 
  (A1) edge (A2)
  (B2) edge (A2)
  (B1) edge (B2)
  (B1) edge (A1)
  (B4) edge (A4)
  (B3) edge (A3)
  (B4) edge (B3)
  (A4) edge (A3);

\node[below] at (A1) {$p$};
\node[below] at (A2) {$q$};
\node[below] at (A3) {$\boldsymbol r$};
\node[below] at (A4) {$q$};
\node[above] at (B1) {$r$};
\node[above] at (B2) {$s$};
\node[above] at (B3) {$t$};
\node[above] at (B4) {$s$};

\end{tikzpicture} &
$  \begin{matrix} \P(p,q^2, r,s^2, t)\\[2mm]   
 \frac1p(1,1),2\times\frac{1}{q}(4,p)\\2\times \frac1s(2,q), \frac1t(1,1)
\end{matrix}$&$ \frac{6 r^3}{r^4-20 r^2+64}        $&$ 1 $ \\

\evnrow $\M\I_{3} $&$ 
  \begin{matrix}
   \boldsymbol {r=2n+9},\\p=r-8,q=r-4\\s=r+4,t=r+8 \\w=(0,4:0,4: p,r)
  \end{matrix}$&
 \begin{tikzpicture}[
  line join=round,
  y={(-0.86cm,0.36cm)},x={(.6cm,0.36cm)}, z={(0cm,.8cm)},
  arr/.style={line cap=round,shorten <= 1pt},baseline={([yshift=-.5em]current bounding box.center)}
]
\def\Side{1}
\coordinate (A1) at (0,0,0);\coordinate (A2) at (0,\Side,0);
\coordinate (A3) at (\Side,\Side,0);\coordinate (A4) at (\Side,0,0);\coordinate (B1) at (0,0,\Side);\coordinate (B2) at (0,\Side,\Side);\coordinate (B3) at (\Side,\Side,\Side);\coordinate (B4) at (\Side,0,\Side);

\draw[thin] (A2) -- (A1) -- (A4);\draw[thin] (B2) -- (B1) -- (B4) -- (B3) -- cycle;\draw[thin] (A1) -- (B1);\draw[thin] (A2) -- (B2);
\draw[thin] (A4) -- (B4);\draw[thin] (A2) -- (A3);
\draw[thin] (A3) -- (B3);\draw[thin] (A3) -- (A4);

\path[arr] 
  (A1) edge (A2)
  (B2) edge (A2)
  (B1) edge (B2)
  (B1) edge (A1)
  (B4) edge (A4)
  (B3) edge (A3)
  (B4) edge (B3)
  (A4) edge (A3);

\node[below] at (A1) {$p$};
\node[below] at (A2) {$r$};
\node[below] at (A3) {$\boldsymbol r$};
\node[below] at (A4) {$q$};
\node[above] at (B1) {$q$};
\node[above] at (B2) {$s$};
\node[above] at (B3) {$t$};
\node[above] at (B4) {$ s$};

\end{tikzpicture} &
$  \begin{matrix} \P(p,q^2, r,s^2, t)\\[2mm]   
 \frac1p(1,1),2\times\frac{1}{q}(8,p)\\2\times \frac1s(4,q), \frac1t(1,1)
\end{matrix}$&$ \frac{6 r^3}{r^4-80 r^2+1024}        $&$ 1$ \\

\oddrow $\M\I_{4} $&$ 
  \begin{matrix}
  \boldsymbol{  r=3n+8},m=2r-1\\p=r-3,q=r-2\\s=r+2,t=2r-3 \\w=(0,2:0,p: 1,r)
  \end{matrix}$&
  \begin{tikzpicture}[
  line join=round,
  y={(-0.86cm,0.36cm)},x={(.6cm,0.36cm)}, z={(0cm,.8cm)},
  arr/.style={line cap=round,shorten <= 1pt},baseline={([yshift=-.5em]current bounding box.center)}
]
\def\Side{1}
\coordinate (A1) at (0,0,0);\coordinate (A2) at (0,\Side,0);
\coordinate (A3) at (\Side,\Side,0);\coordinate (A4) at (\Side,0,0);\coordinate (B1) at (0,0,\Side);\coordinate (B2) at (0,\Side,\Side);\coordinate (B3) at (\Side,\Side,\Side);\coordinate (B4) at (\Side,0,\Side);

\draw[thin] (A2) -- (A1) -- (A4);\draw[thin] (B2) -- (B1) -- (B4) -- (B3) -- cycle;\draw[thin] (A1) -- (B1);\draw[thin] (A2) -- (B2);
\draw[thin] (A4) -- (B4);\draw[thin] (A2) -- (A3);
\draw[thin] (A3) -- (B3);\draw[thin] (A3) -- (A4);

\path[arr] 
  (A1) edge (A2)
  (B2) edge (A2)
  (B1) edge (B2)
  (B1) edge (A1)
  (B4) edge (A4)
  (B3) edge (A3)
  (B4) edge (B3)
  (A4) edge (A3);

\node[below] at (A1) {$1$};
\node[below] at (A2) {$q$};
\node[below] at (A3) {$\boldsymbol r$};
\node[below] at (A4) {$3$};
\node[above] at (B1) {$r$};
\node[above] at (B2) {$t$};
\node[above] at (B3) {$m$};
\node[above] at (B4) {$s$};

\end{tikzpicture} &
 $\begin{matrix} \P(1,3,q, r,s, t,m)\\[2mm]   
 \frac13(1,1),\dfrac{1}{q}(1,1)\\ \frac1s(3,p), \frac1t(2,q),\frac1m(s,t)
\end{matrix}$&$ \frac{2 r^3 \left(2 r^2+4 r-7\right)}{3 (r-2) (r+2) (2 r-3) (2 r-1)}        $&$ 4$ \\

\evnrow  $\M\I_{5} $&$ 
  \begin{matrix}
  \boldsymbol{  r=6n+2},\\p=r-2,q=r-1\\s=r+1,t=r+2 \\w=(0,1:0,1: p,r)
  \end{matrix}$&
 \begin{tikzpicture}[
  line join=round,
  y={(-0.86cm,0.36cm)},x={(.6cm,0.36cm)}, z={(0cm,.8cm)},
  arr/.style={line cap=round,shorten <= 1pt},baseline={([yshift=-.5em]current bounding box.center)}
]
\def\Side{1}
\coordinate (A1) at (0,0,0);\coordinate (A2) at (0,\Side,0);
\coordinate (A3) at (\Side,\Side,0);\coordinate (A4) at (\Side,0,0);\coordinate (B1) at (0,0,\Side);\coordinate (B2) at (0,\Side,\Side);\coordinate (B3) at (\Side,\Side,\Side);\coordinate (B4) at (\Side,0,\Side);

\draw[thin] (A2) -- (A1) -- (A4);\draw[thin] (B2) -- (B1) -- (B4) -- (B3) -- cycle;\draw[thin] (A1) -- (B1);\draw[thin] (A2) -- (B2);
\draw[thin] (A4) -- (B4);\draw[thin] (A2) -- (A3);
\draw[thin] (A3) -- (B3);\draw[thin] (A3) -- (A4);

\path[arr] 
  (A1) edge (A2)
  (B2) edge (A2)
  (B1) edge (B2)
  (B1) edge (A1)
  (B4) edge (A4)
  (B3) edge (A3)
  (B4) edge (B3)
  (A4) edge (A3);

\node[below] at (A1) {$p$};
\node[below] at (A2) {$q$};
\node[below] at (A3) {$\boldsymbol r$};
\node[below] at (A4) {$q$};
\node[above] at (B1) {$r$};
\node[above] at (B2) {$s$};
\node[above] at (B3) {$t$};
\node[above] at (B4) {$ s$};

\end{tikzpicture} &
  $\begin{matrix} \P(p,q^2, r,s^2, t)\\[2mm]   
 \frac1p(1,1),2\times\frac{1}{q}(2,p)\\2\times \frac1s(1,q), \frac1t(1,1)
\end{matrix}$&$ \frac{6 r^3}{r^4-5 r^2+4}       $&$ 1$ \\
\end{longtable}\end{thm}
\begin{thm}\label{Th:Weights}\textbf{(Weight models)}  There exist at least seven  biregular  weight models  of wellformed and quasismooth rigid  del Pezzo surfaces of Fano index 1 or 2   such that  their images under the anti-canonical  or sub anti-canonical embeddings in   \(\P^6(a,\ldots,g)\) are modeled on those of the  Segre embedding of \(\PI \) in \(\P^7\).  At least, five of them do not deform under  $\Q$-Gorenstein degeneration to a toric variety. \label{thrm-pxp1} 

 \renewcommand*{\arraystretch}{1.2}
\begin{longtable}{>{\hspace{0.5em}}llccccr<{\hspace{0.5em}}}

\caption{biregular weight models of rigid del  Pezzo surfaces of index 1 and 2 in  $\PI $ format where $p=r-2, q=r-1, s=r+1,t=r+2, m=2r-1, u=2r, v=2r+1,z=3r-1$. The \(j\)-th model of the Fano index \(i\) is denoted by \(\M_{ij}\).} \label{Weight-Models}\\
\toprule
\multicolumn{1}{c}{Model}&\multicolumn{1}{c}{WPS \& Para\ }&\multicolumn{1}{c}{Weight Cube}&\multicolumn{1}{c}{$\mathcal{B}$}&\multicolumn{1}{c}{$-K^2$}&\multicolumn{1}{c}{$h^0(-K)$}\\
\cmidrule(lr){1-1}\cmidrule(lr){2-2}\cmidrule(lr){3-3}\cmidrule(lr){4-4}\cmidrule(lr){5-6}
\endfirsthead
\multicolumn{7}{l}{\vspace{-0.25em}\scriptsize\emph{\tablename\ \thetable{} continued from previous page}}\\
\midrule
\endhead
\multicolumn{7}{r}{\scriptsize\emph{Continued on next page}}\\
\endfoot
\bottomrule
\endlastfoot

\evnrow $\M_{11} $&$ 
  \begin{matrix}
  \P(1,2, r^2,s^2, m)\\  \; r=2n+2, n\ge 1\\w=(0,1:0,q: 1,r)
  \end{matrix}$&$
 \begin{tikzpicture}[
  line join=round,
  y={(-0.86cm,0.36cm)},x={(.6cm,0.36cm)}, z={(0cm,.8cm)},
  arr/.style={line cap=round,shorten <= 1pt},baseline={([yshift=-.5em]current bounding box.center)}
]
\def\Side{1}
\coordinate (A1) at (0,0,0);\coordinate (A2) at (0,\Side,0);
\coordinate (A3) at (\Side,\Side,0);\coordinate (A4) at (\Side,0,0);\coordinate (B1) at (0,0,\Side);\coordinate (B2) at (0,\Side,\Side);\coordinate (B3) at (\Side,\Side,\Side);\coordinate (B4) at (\Side,0,\Side);

\draw[thin] (A2) -- (A1) -- (A4);\draw[thin] (B2) -- (B1) -- (B4) -- (B3) -- cycle;\draw[thin] (A1) -- (B1);\draw[thin] (A2) -- (B2);
\draw[thin] (A4) -- (B4);\draw[thin] (A2) -- (A3);
\draw[thin] (A3) -- (B3);\draw[thin] (A3) -- (A4);

\path[arr] 
  (A1) edge (A2)
  (B2) edge (A2)
  (B1) edge (B2)
  (B1) edge (A1)
  (B4) edge (A4)
  (B3) edge (A3)
  (B4) edge (B3)
  (A4) edge (A3);

\node[below] at (A1) {$1$};
\node[below] at (A2) {$r$};
\node[below] at (A3) {$s$};
\node[below] at (A4) {$2$};
\node[above] at (B1) {$r$};
\node[above] at (B2) {$m$};
\node[above] at (B3) {$\boldsymbol {{u}}$};
\node[above] at (B4) {$s$};

\end{tikzpicture}$ &
  $    
2\times \dfrac1s(2,r),\dfrac{1}{m}(1,1)
$&$ \dfrac{2 r+5}{2 r^2+r-1}         $&$ 1$ \\

\oddrow $\M_{12} $&$ 
  \begin{matrix}
  \P(2, r,s^{3},u, z)\\  \; r=2n,\; n\ge 1\\
  w=(0,q:0,q: 2,s)
  \end{matrix}$&
 $\begin{tikzpicture}[
  line join=round,
  y={(-0.86cm,0.36cm)},x={(.6cm,0.36cm)}, z={(0cm,.8cm)},
  arr/.style={line cap=round,shorten <= 1pt},baseline={([yshift=-.5em]current bounding box.center)}
]
\def\Side{1}
\coordinate (A1) at (0,0,0);\coordinate (A2) at (0,\Side,0);
\coordinate (A3) at (\Side,\Side,0);\coordinate (A4) at (\Side,0,0);\coordinate (B1) at (0,0,\Side);\coordinate (B2) at (0,\Side,\Side);\coordinate (B3) at (\Side,\Side,\Side);\coordinate (B4) at (\Side,0,\Side);

\draw[thin] (A2) -- (A1) -- (A4);\draw[thin] (B2) -- (B1) -- (B4) -- (B3) -- cycle;\draw[thin] (A1) -- (B1);\draw[thin] (A2) -- (B2);
\draw[thin] (A4) -- (B4);\draw[thin] (A2) -- (A3);
\draw[thin] (A3) -- (B3);\draw[thin] (A3) -- (A4);

\path[arr] 
  (A1) edge (A2)
  (B2) edge (A2)
  (B1) edge (B2)
  (B1) edge (A1)
  (B4) edge (A4)
  (B3) edge (A3)
  (B4) edge (B3)
  (A4) edge (A3);

\node[below] at (A1) {$2$};
\node[below] at (A2) {$s$};
\node[below] at (A3) {$\boldsymbol u$};
\node[below] at (A4) {$s$};
\node[above] at (B1) {$s$};
\node[above] at (B2) {$u$};
\node[above] at (B3) {$z$};
\node[above] at (B4) {$\boldsymbol u$};
\end{tikzpicture} $ &$
3\times\dfrac{1}{s}(2,r),\dfrac{1}{z}(r,u)
$&$\dfrac{6}{(r+1) (3 r-1)}
 $&$ 0$ \\

\evnrow $\M_{13} $&$ 
  \begin{matrix}
  \P(2,3, r^2,s,t, m)\\  \; r=6n-1, n\ge 1\\w=(0,1:0,p: 2,s)
  \end{matrix}$&
 \begin{tikzpicture}[
  line join=round,
  y={(-0.86cm,0.36cm)},x={(.6cm,0.36cm)}, z={(0cm,.8cm)},
  arr/.style={line cap=round,shorten <= 1pt},baseline={([yshift=-.5em]current bounding box.center)}
]
\def\Side{1}
\coordinate (A1) at (0,0,0);\coordinate (A2) at (0,\Side,0);
\coordinate (A3) at (\Side,\Side,0);\coordinate (A4) at (\Side,0,0);\coordinate (B1) at (0,0,\Side);\coordinate (B2) at (0,\Side,\Side);\coordinate (B3) at (\Side,\Side,\Side);\coordinate (B4) at (\Side,0,\Side);

\draw[thin] (A2) -- (A1) -- (A4);\draw[thin] (B2) -- (B1) -- (B4) -- (B3) -- cycle;\draw[thin] (A1) -- (B1);\draw[thin] (A2) -- (B2);
\draw[thin] (A4) -- (B4);\draw[thin] (A2) -- (A3);
\draw[thin] (A3) -- (B3);\draw[thin] (A3) -- (A4);

\path[arr] 
  (A1) edge (A2)
  (B2) edge (A2)
  (B1) edge (B2)
  (B1) edge (A1)
  (B4) edge (A4)
  (B3) edge (A3)
  (B4) edge (B3)
  (A4) edge (A3);

\node[below] at (A1) {$2$};
\node[below] at (A2) {$s$};
\node[below] at (A3) {$t$};
\node[below] at (A4) {$3$};
\node[above] at (B1) {$r$};
\node[above] at (B2) {$m$};
\node[above] at (B3) {$\boldsymbol u$};
\node[above] at (B4) {$\boldsymbol s$};

\end{tikzpicture} &
  $    
\begin{matrix}\frac13(1,1),3\times \frac1r(2,q),\\\frac1t(3,r),\frac{1}{m}(r,r)\end{matrix}
$&$ \dfrac{(r+1) (2 r+9)}{3 r (r+2) (2 r-1)}        $&$ 0$ \\
\oddrow $\M_{14} $&$ 
  \begin{matrix}
  \P(2,3, r^2,s,t, m)\\  \; r=6n+3, n\ge 1\\w=(0,1:0,p: 2,s)
  \end{matrix}$&\begin{tikzpicture}[
  line join=round,
  y={(-0.86cm,0.36cm)},x={(.6cm,0.36cm)}, z={(0cm,.8cm)},
  arr/.style={line cap=round,shorten <= 1pt},baseline={([yshift=-.5em]current bounding box.center)}
]
\def\Side{1}
\coordinate (A1) at (0,0,0);\coordinate (A2) at (0,\Side,0);
\coordinate (A3) at (\Side,\Side,0);\coordinate (A4) at (\Side,0,0);\coordinate (B1) at (0,0,\Side);\coordinate (B2) at (0,\Side,\Side);\coordinate (B3) at (\Side,\Side,\Side);\coordinate (B4) at (\Side,0,\Side);

\draw[thin] (A2) -- (A1) -- (A4);\draw[thin] (B2) -- (B1) -- (B4) -- (B3) -- cycle;\draw[thin] (A1) -- (B1);\draw[thin] (A2) -- (B2);
\draw[thin] (A4) -- (B4);\draw[thin] (A2) -- (A3);
\draw[thin] (A3) -- (B3);\draw[thin] (A3) -- (A4);

\path[arr] 
  (A1) edge (A2)
  (B2) edge (A2)
  (B1) edge (B2)
  (B1) edge (A1)
  (B4) edge (A4)
  (B3) edge (A3)
  (B4) edge (B3)
  (A4) edge (A3);

\node[below] at (A1) {$2$};
\node[below] at (A2) {$s$};
\node[below] at (A3) {$t$};
\node[below] at (A4) {$3$};
\node[above] at (B1) {$r$};
\node[above] at (B2) {$m$};
\node[above] at (B3) {$\boldsymbol u$};
\node[above] at (B4) {$\boldsymbol s$};

\end{tikzpicture} &
  $    
\begin{matrix}2\times \frac13(1,1),3\times \frac1r(2,q),\\\frac1t(3,r),\frac{1}{m}(r,r)\end{matrix}
$&same as \(\M_{13}\)&$ 0$ \\
 \evnrow $\M_{15} $&$ 
  \begin{matrix}
  \P(1,3, r^2,t^{2}, m)\\  \; r=6n+3, n\ge 1\\w=(0,2:0,q: 1,r)
  \end{matrix}$&
  \begin{tikzpicture}[
  line join=round,
  y={(-0.86cm,0.36cm)},x={(.6cm,0.36cm)}, z={(0cm,.8cm)},
  arr/.style={line cap=round,shorten <= 1pt},baseline={([yshift=-.5em]current bounding box.center)}
]
\def\Side{1}
\coordinate (A1) at (0,0,0);\coordinate (A2) at (0,\Side,0);
\coordinate (A3) at (\Side,\Side,0);\coordinate (A4) at (\Side,0,0);\coordinate (B1) at (0,0,\Side);\coordinate (B2) at (0,\Side,\Side);\coordinate (B3) at (\Side,\Side,\Side);\coordinate (B4) at (\Side,0,\Side);

\draw[thin] (A2) -- (A1) -- (A4);\draw[thin] (B2) -- (B1) -- (B4) -- (B3) -- cycle;\draw[thin] (A1) -- (B1);\draw[thin] (A2) -- (B2);
\draw[thin] (A4) -- (B4);\draw[thin] (A2) -- (A3);
\draw[thin] (A3) -- (B3);\draw[thin] (A3) -- (A4);

\path[arr] 
  (A1) edge (A2)
  (B2) edge (A2)
  (B1) edge (B2)
  (B1) edge (A1)
  (B4) edge (A4)
  (B3) edge (A3)
  (B4) edge (B3)
  (A4) edge (A3);

\node[below] at (A1) {$1$};
\node[below] at (A2) {$r$};
\node[below] at (A3) {$t$};
\node[below] at (A4) {$3$};
\node[above] at (B1) {$r$};
\node[above] at (B2) {$m$};
\node[above] at (B3) {$\boldsymbol v$};
\node[above] at (B4) {$t$};

\end{tikzpicture} &
  $    
\begin{matrix}\frac13(1,1),2\times \frac1r(2,q),\\2\times\frac1t(3,q),\frac{1}{m}(r,r)\end{matrix}
$&$ \dfrac{2 \left(2 r^2+8 r+3\right)}{3 r (r+2) (2 r-1)}       $&$ 1$ \\\hline\hline
\oddrow $\M_{21} $&$ 
  \begin{matrix}
  \P(3, q,r^2,s^3)\\  \; r=3n+3, n\ge 1\\w=(0,1:0,1: q,r)
  \end{matrix}$&
  \begin{tikzpicture}[
  line join=round,
  y={(-0.86cm,0.36cm)},x={(.6cm,0.36cm)}, z={(0cm,.8cm)},
  arr/.style={line cap=round,shorten <= 1pt},baseline={([yshift=-.5em]current bounding box.center)}
]
\def\Side{1}
\coordinate (A1) at (0,0,0);\coordinate (A2) at (0,\Side,0);
\coordinate (A3) at (\Side,\Side,0);\coordinate (A4) at (\Side,0,0);\coordinate (B1) at (0,0,\Side);\coordinate (B2) at (0,\Side,\Side);\coordinate (B3) at (\Side,\Side,\Side);\coordinate (B4) at (\Side,0,\Side);

\draw[thin] (A2) -- (A1) -- (A4);\draw[thin] (B2) -- (B1) -- (B4) -- (B3) -- cycle;\draw[thin] (A1) -- (B1);\draw[thin] (A2) -- (B2);
\draw[thin] (A4) -- (B4);\draw[thin] (A2) -- (A3);
\draw[thin] (A3) -- (B3);\draw[thin] (A3) -- (A4);

\path[arr] 
  (A1) edge (A2)
  (B2) edge (A2)
  (B1) edge (B2)
  (B1) edge (A1)
  (B4) edge (A4)
  (B3) edge (A3)
  (B4) edge (B3)
  (A4) edge (A3);

\node[below] at (A1) {$q$};
\node[below] at (A2) {$r$};
\node[below] at (A3) {$s$};
\node[below] at (A4) {$\boldsymbol r$};
\node[above] at (B1) {$r$};
\node[above] at (B2) {$s$};
\node[above] at (B3) {$\boldsymbol t$};
\node[above] at (B4) {$s$};

\end{tikzpicture} &
  $    
\begin{matrix}3\times \frac13(1,1), \frac1q(1,1),\\3\times \frac1s(3,r)\end{matrix}
$&$ \dfrac{8}{r^2-1}        $&$ 0$ \\
\evnrow $\M_{22} $&$ 
  \begin{matrix}
  \P(3, p,q^{2},r,s,t)\\  \; r=3n+5, n\ge 1\\w=(0,1:0,1: p,r)
  \end{matrix}$&
 \begin{tikzpicture}[
  line join=round,
  y={(-0.86cm,0.36cm)},x={(.6cm,0.36cm)}, z={(0cm,.8cm)},
  arr/.style={line cap=round,shorten <= 1pt},baseline={([yshift=-.5em]current bounding box.center)}
]
\def\Side{1}
\coordinate (A1) at (0,0,0);\coordinate (A2) at (0,\Side,0);
\coordinate (A3) at (\Side,\Side,0);\coordinate (A4) at (\Side,0,0);\coordinate (B1) at (0,0,\Side);\coordinate (B2) at (0,\Side,\Side);\coordinate (B3) at (\Side,\Side,\Side);\coordinate (B4) at (\Side,0,\Side);

\draw[thin] (A2) -- (A1) -- (A4);\draw[thin] (B2) -- (B1) -- (B4) -- (B3) -- cycle;\draw[thin] (A1) -- (B1);\draw[thin] (A2) -- (B2);
\draw[thin] (A4) -- (B4);\draw[thin] (A2) -- (A3);
\draw[thin] (A3) -- (B3);\draw[thin] (A3) -- (A4);

\path[arr] 
  (A1) edge (A2)
  (B2) edge (A2)
  (B1) edge (B2)
  (B1) edge (A1)
  (B4) edge (A4)
  (B3) edge (A3)
  (B4) edge (B3)
  (A4) edge (A3);

\node[below] at (A1) {$p$};
\node[below] at (A2) {$r$};
\node[below] at (A3) {$\boldsymbol s$};
\node[below] at (A4) {$q$};
\node[above] at (B1) {$q$};
\node[above] at (B2) {$s$};
\node[above] at (B3) {$t$};
\node[above] at (B4) {$\boldsymbol r$};

\end{tikzpicture} &
 $     
\begin{matrix}3\times \frac13(1,1), \frac1p(1,1),\\\frac1q(3,p), \frac1t(3,s)\end{matrix}$
&$ \dfrac{8r}{(r^2-4) (r-1) }     $&$ 0$ 
\end{longtable}
 \end{thm}

This paragraph explains the notations appearing in the Table in  \ref{Index-Models}. and \ref{Weight-Models}. All of our families of rigid odPs can be considered as weighted complete intersections of some weighted \(\PI\) variety (defined in section \ref{S-P1^3}) or of some (weighted) projective cone(s) over it.   The notation WPS gives the weights of ambient  weighted projective space  \(\P^6(w_i)\) containing  our odPs, with powers of \(w_{i}\) representing the multiplicity of the given weight.   The weights  on \(\P^7(w_i)  \) containing  a w\((\PI)\) variety are determined by a choice of parameter \(w\). The cube in the given table represents the weights of the ambient weighted \( \PI\) variety.  The bold numbers in the weight cube represent the degree(s) of  quasilinear section(s) of \( \P^{7}(w_i)\) with whom we take the intersection. If  some  numbers \(a_is \) appear in \(\P^6(w_i)\) in the column WPS which are not in  the weight cube, then our surface appears as a quasilinear section of some (weighted) projective cone(s) over given w\((\PI)\) variety, with cone(s) having weights \(a_is\).  The columns \(\mathcal B, -K^2 \) and \(h^0(-K)\) lists the basket of rigid singularities,  the anti-canonical degrees, and the first plurigenus  of \(X \), respectively. 

The first key  point of our computation is to prove the general Hilbert series formula for the a weighted \( \PI\) variety, which enables us to compute the anti-canonical divisor class of the ambient variety. Then we use an  algorithmic method \cite{QJSC,BKZ}  to search for candidate rigid del Pezzo surfaces by using the computer algebra system \textsc{MAGMA} ~\cite{magma}.   We analyze these candidates by comparing the basket of orbifold points
and weights
  to spot any pattern emerging among them which may lead to a biregular model (infinite series)  of such  surfaces. Once such a model is identified,  we use the strategy described in section \ref{S-Proof} to prove the existence,  wellformedness and quasismoothness of the given models. The list of biregular models  is non-exhaustive due to the non existence of termination conditions on the algorithmic approach used to compute candidate del Pezzo surfaces. However, our list is the complete classification  up to certain computational bounds described in section \ref{computer-search}. In the Fano index 1 case, there are also wellformed and quasismooth families of rigid odPs that  are not part of any model, and they are listed as sporadic cases in section \ref{Sporadic}.

\subsection{Context and motivation}
The classification of non-singular del Pezzo surfaces is well understood, and there are 10 deformation families. An orbifold del Pezzo surface \(X\), classically known as a log del Pezzo surface, is a del Pezzo surface with at worst isolated orbifold points. If \(X\) contains only \(R\) (rigid) singularities, then it is called \(\Q\)-Gorenstein (\(qG\))-rigid. It is said to be of class TG if it can be deformed under a \(qG\)-degeneration to a normal toric del Pezzo surface. Since all normal surfaces with orbifold points have partial smoothing to a surface containing only \(R\) singularities \cite{Kollar-SB}, our focus in this article is only on rigid del Pezzo surfaces.

Orbifold del Pezzo surfaces appear in various veins in algebraic geometry, for example, classically in the minimal model program~\cite{Kawamata-MMP}. More recently, the classification of orbifold del Pezzo surfaces gained traction due to the mirror symmetry program \cite{CCGGK} of Coates--Corti--Galkin--Kaspezyk. In \cite{dp-AMS}, it is conjectured that there is a one-to-one correspondence between mutation equivalent classes of Fano polygons and \(qG\)-deformation equivalence classes of rigid orbifold del Pezzo surfaces. This conjecture inspired the classification of orbifold del Pezzo surfaces with \(k \times \frac{1}{3} (1,1)\) singularities by Corti--Heuberger \cite{dp-CH}. The article provided explicit constructions of these varieties as complete intersections of so-called rep-quotient varieties (toric varieties, weighted Grassmannians, etc.), to be able to compute the quantum orbifold cohomology and give evidence for the conjecture. The explicit constructions of such surfaces in specific Gorenstein formats are a particular instance of finding complete intersection models of del Pezzo surfaces in rep-quotient varieties, directly linking to mirror symmetry.

The classification with a single orbifold point of type \(\frac{1}{r}(1,1)\) has been given by Cavey--Prince \cite{dp-CP}, and the classification with a baskets of orbifold points of type \[\cB = \left\{ k_1 \times \frac{1}{3}(1,1) + k_2 \times \frac{1}{5}(1,2) \right\} \text{ and  } h^0(-K_X) \ne 0\] has been settled in the Ph.D. thesis of Cuzzucoli \cite{Cuzzucoli}. The thesis lists possible candidate surfaces, having a baskets  of type \(\cB\)  and  \(h^0(-K_X) = 0\). We prove the existence of two these  of rigid del Pezzo surfaces in \(\PI\) format: number 5 in Table \ref{Tab:Sporadic} \((k_1 = 0, k_2 = 3)\) and the case \(n = 1\) in the model \(\M_{12}\) \((k_1 = 3, k_2 = 1)\).

Another context where the explicit descriptions of orbifold del Pezzo surfaces embedded in some weighted projective space have been used is in the study of \(K\)-stability questions that require the computation of the log-canonical threshold to examine the existence of a K\"ahler-Einstein metric on such varieties. This was studied by Johnson--Koll\'ar \cite{Jhonson-Kollar} for certain classes of del Pezzo hypersurfaces and further investigated for other classes of hypersurfaces in \cite{dp-exceptional, dp-zoo} and for complete intersections in \cite{Park}. We hope that the explicit non-complete intersection descriptions of orbifold del Pezzo surfaces that we provide can be utilized in this vein.

Our constructions can also be seen as an extension of existing classification results of describing orbifold del Pezzo surfaces as low codimension varieties
in some weighted projective space. The complete classification of well-formed and quasi-smooth del Pezzo hypersurfaces was provided by Paemurru \cite{dp-Erik}
and  combinatorial classification of all ordered tuples \((d_1, d_2, w_1, w_2, w_3, w_4, w_5)\) such that \(X_{d_1,d_2} \into \P(w_1, \ldots, w_5)\) is a well-formed and quasi-smooth weighted complete intersection del Pezzo surface was classified by Mayanskiy \cite{Mayanskiy}. The classification of low codimension del Pezzo surfaces with specific bounds on the orbifold points appeared in \cite{Q-Mathematics}. There are other results on the classification of log del Pezzo surfaces by Nikulin \cite{Niku1, Niku2, Niku3}, Nakayama \cite{Nakayama}, and Hwang--Park \cite{Hwang-Park}. The construction of biregular models of orbifold del Pezzo surfaces in the non-complete intersection case was initiated in \cite{QMOC}, which we extend further in this paper.

Last but not least, our main emphasis in this paper is on the constructions
of  orbifold del Pezzo surfaces. However, weighted \(\PI\) is one of the most well-understood codimension 4 Gorenstein variety, and the Hilbert series formula \eqref{HS-P1^3} and corresponding graded ring structure can be used to construct other interesting classes of polarized orbifolds, such as Calabi--Yau 3-folds, \(\Q\)-Fano 3-folds, and 3-folds of general type.

\subsection{Summary and Notation}
 We describe the notion of a weighted \(\PI\) variety in section \ref{S-P1^3} and give a general formula for its Hilbert Series. The  section \ref{S-Proof} provides detailed steps of our proofs and explains various intricacies appearing during our calculations. We also provide a detailed explanation in two sample biregular models and further explanation of our calculation in   section \ref{S-Proof}.
The  section \ref{S-Sporadic} gives details of our computer search and lists all the sporadic cases of rigid orbifold del Pezzo surfaces constructed in this paper.

We use the following notation to denote various mathematical objects in
the rest of the article.

\begin{itemize}

\item For a simplicity of notation, we write  
$\frac 1r(1,b')$  for an orbifold point type \(\frac 1r(a,b)\) in the tables
where we specify the orbifold points.  This is always possible by choosing
 a different
generator of the cyclic group \(\mu_r\), and was proved in \cite{QMOC}. 
\item  We use   $\wP$
to  denote the ambient weighted $\PI$ variety. 
\item $a,b,c,d,e,f$ and $g$ are used to represent the variables of the ambient space \(\P^6(w_i)\), and the weights are represented by  $m,r,s,t,u,v,y, \textrm{ and }z$.  The notation \(a_m\) will mean that the variable \(a\) has weight \(m\).   
\item We use capitals $H_ds \text{ and } J_ds$ to denote the weighted homogeneous forms appearing in the equations of the del Pezzo surface \(X\). 
\end{itemize}  

%--------------------------------------------------------------------------------------------
 %%%%%%%%%%%%%%%%

%--------------------------------------------------------------------------------------------

 %%%%%%%%%%%%%%%%%%%%%%%%%%%%%%%%%%%%
\section{ Weighted $\mathbb P^1 \times \mathbb P^1 \times\mathbb P^1  $   }
\label{S-P1^3}

In this section, we describe the  construction of weighted \(\PI\) variety \(\wP\) and give a general formula for its Hilbert series. 
The formula allows us to compute the graded free resolution and orbifold canonical class of \(\wP.\)   Consider the  Segre embedding of $\PI$   
\[\P_{s_1,s_2}\times \P_{t_1,t_2} \times \P_{u_1,u_2} \into \P^7\left(s_it_ju_k\right)=\P\left(X_{ijk}\right), 1\le i,j,k\le 2 .\] Then the image of  $\PI\into \P^7$ can be described  by the $2 \times 2$ minors of the cube 

 \begin{equation}\label{eq!cube}\begin{tikzpicture}[
  line join=round,
  y={(-0.86cm,0.36cm)},x={(.6cm,0.36cm)}, z={(0cm,.8cm)},
  arr/.style={line cap=round,shorten <= 1pt},baseline={([yshift=-.5em]current bounding box.center)}
]
\def\Side{2}
\coordinate (A1) at (0,0,0);\coordinate (A2) at (0,\Side,0);
\coordinate (A3) at (\Side,\Side,0);\coordinate (A4) at (\Side,0,0);\coordinate (B1) at (0,0,\Side);\coordinate (B2) at (0,\Side,\Side);\coordinate (B3) at (\Side,\Side,\Side);\coordinate (B4) at (\Side,0,\Side);
\draw[thin] (A2) -- (A1) -- (A4);\draw[thin] (B2) -- (B1) -- (B4) -- (B3) -- cycle;\draw[thin] (A1) -- (B1);\draw[thin] (A2) -- (B2);
\draw[thin] (A4) -- (B4);\draw[thin] (A2) -- (A3);
\draw[thin] (A3) -- (B3);\draw[thin] (A3) -- (A4);
\path[arr] 
  (A1) edge (A2)
  (B2) edge (A2)
  (B1) edge (B2)
  (B1) edge (A1)
  (B4) edge (A4)
  (B3) edge (A3)
  (B4) edge (B3)
  (A4) edge (A3);

\node[below] at (A1) {\;\;\; $X_{111}$};
\node[below] at (A2) {\hspace{-5mm}$X_{121}$};
\node[below] at (A3) {\;\;\;$X_{221}$};
\node[below] at (A4) {\;\;\;\;\;$X_{211}$};
\node[above] at (B1) {\hspace{-2mm}$X_{112}$};
\node[above] at (B2) {\hspace{-6mm}$X_{122}$};
\node[above] at (B3) {$X_{222}$};
\node[above] at (B4) {\hspace{6mm}$X_{212}$};
\end{tikzpicture}.
\end{equation}

It is well known that it is  a Gorenstein variety of codimension 4 having a \(9\times 16\) resolution \cite{BKR}.    

%\subsection{Weighted $\PxP$}
%\label{s!sPxP}
We use the structure of $\P^1$ as a quotient of $\GL(2,\C)$ by its group of upper triangular matrices $B$ of full rank to define  the weighted $\PI$ variety, following Corti--Reid~\cite{wg} and  \cite{QS,QS-AHEP,BKQ}.

Let  $G = \GL(2,\C)\times \GL(2,\C) \times \GL(2,\C)$ and $B=B_1\times B_2 \times B_3$ is the  product of the corresponding groups of upper triangular matrices. Then $G$ has rank 6 character lattice  $\Lambda=\Hom(T,\C^*)\cong\Z^6$, for the maximal torus $T\subset G$.
We choose a vector $$\mu=(a_1,a_2,b_1   ,b_2,c_1,c_2)\in \Hom(\Lambda,\Z),$$   such that $a_i+b_j+c_k>0$ for \(1\le i,j,k\le2\).
Then we define weighted $\PI$ by  \begin{equation}
\wP=\left(C(\PI)\setminus\{0\}\right)/\C^\times
\end{equation}  where $\C^\times$ acts on $C(\PI)\subset \C^8$ by $$X_{ijk}\mapsto \lambda^{a_i+b_j+c_k}X_{ijk},\;\; i,j,k \in \{1,2\}. $$ Thus we have the embedding
\begin{equation}\label{e!wSigma}
\wP \hookrightarrow w\P^7\left(a_i+b_j+c_k\right), i,j,k\in \{1,2\} 
\end{equation}
with image defined by $2\times 2$ minors of the weight  cube \eqref{eq!cube}  with respect to the weights 
\begin{equation}
\deg\left(X_{ikj}\right)=a_i+b_j+c_k,\;\; \forall i,j,k. \end{equation}
By using the  general Hilbert series
formula of~\cite[Theorem 3.1]{QS}, we get the following formula for the Hilbert series of general \(\wP\).

\begin{prop}
The Hilbert series of $\wP$ in the embedding  \eqref{e!wSigma} is 
\begin{equation}
\label{HS-P1^3}
P_{\wP}(t) =
\dfrac{1-Q_1t^s+Q_2t^s-Q_1t^{2s}+t^{3s}}{\displaystyle\prod_{1\le i,j,k\le 2}\left(1-t^{a_i+b_j+c_k}\right)},
\end{equation}
where 
\[Q_1=3+\sum_{1 \le i\ne j \le 2} t^{a_i-a_j}+\sum_{1 \le i\ne j \le 2} t^{b_i-b_j}+\sum_{1 \le i\ne j \le 2} t^{c_i-c_j},\;\;\; Q_2=\sum_{1\le i,j,k \le 2}2t^{a_i+b_j+c_k},
\]
 and  $l=a_1+a_2+b_1+b_2+c_1+c_2$. Moreover, if $\wP$ is wellformed then $$K_{\wP}=\cO(-l).$$
\end{prop}
The first two terms in the numerator represent the  $9\times 16$ resolution.
For instance, the term $t^{a_1-a_2+s} = t^{(a_1+b_1+c_1) + (a_1+b_2+c_2)}$  has the degree of the equation $X_{111}X_{122}-X_{112}X_{121}$.

\begin{eg} For $\mu=(0,0:0,0:1,1)$  we get the straight  $\PI\into\P^7$ with the canonical bundle $\cO(-2)$ which is a del Pezzo 3-fold of index 2. The hyperplane section of it gives a family describing  classical smooth del Pezzo surface of degree 6.  
\end{eg}
\begin{eg} For $\mu=(0,1:0,1:1,2)$  we get the embedding   $$\wP\into \P^6(1,2^3,3^3,4) $$ with the canonical bundle $\cO(-5)$. The intersection with a generic quartic \[X:=\wP\cap (4)\into \P(1,2^3,3^3)\]  is a family of orbifold del Pezzo surface of index 1 with 3 orbifold points of type  \(\frac 13(1,1)\). The toric construction for this surface was given in \cite{dp-CH}.   
\end{eg}

%%%%%%%%%%%%%%%%%%%%%%%%%%%%%%%%%%%%%

\section{Proofs of the  results}
\label{S-Proof}
We first  describe the general strategy of proving the main theorem and  list all the steps required in the proof of the main results. Then we describe the proof of main theorems in representative cases. 
\subsection{General strategy of proofs}
\label{proof-strategy}
For each model, and the sporadic examples of section \ref{Sporadic}, the proofs are divided into several steps. We briefly discuss them here as the detailed explanation  can be found in \cite{QMOC}.  
 \subsubsection{Computer search algorithm }
 \label{algorithm}
 This section recalls the algorithm from \cite{QJSC} which we used to compute the list of candidate rigid del Pezzo surfaces.  An important ingredient is the  Riemann--Roch formula
 \cite{BRZ} which gives a  formula for Hilbert series \(P_X(t)\) of \(X\) as the sum of  a smooth part and an orbifold part.
It  states that if   an orbifold \(X\)  has a collection \(\mathcal B=\{k_i \times Q_i: m_i \in \Z_{>0}\}\) of isolated orbifold points then;
\begin{equation}\label{hs-decomp}P_X(t)=P_\sm(t)+\sum_{Q_i \in \mathcal B} k_iP_{Q_i}(t),\end{equation}
where \(\P_{\sm}\) gives the smooth part and represents the orbifold part of the Hilbert Series. Indeed if the variety is smooth then \(P_X(t)=P_{\sm }(t)\). The algorithm provides a complete list of  orbifolds having  fixed Fano index, dimension, and  canonical class \( K_{X}=\Oh(kD)\) in a given ambient  space. Indeed, if  \(X\) is a del Pezzo surface  of index \(I\) then \(k=-I. \) We briefly recall the steps  of our algorithm for the case of del Pezzo surfaces. 
\begin{enumerate}
\item Calculate  the Hilbert series and the canonical divisor class of the of  \(\wP\).
\item 
Find all possible embeddings  in \(\P^6(w_i)\) of index \(I\), i.e.   \(K_X=\Oh(-I)\) by  enumerating weights and using the adjunction formula.
\item 
For each  embedding computed at step (ii), calculate the Hilbert series of \(X\)  and  the term representing its  smooth part \(P_\sm(t)\). 
\item Compute the list of all possible   rigid orbifold points from the weights of the \(\P^6(w_i)\) that may lie on \(X\).
\item Form all possible subsets of orbifold points and for each set \(\mathcal B\),   determine if the difference  \(P_X(t)-P_{\sm}(t)\) can be equal to  $\sum_{Q_i \in \mathcal B} k_iP_{Q_i}(t)$  for some \(k_{i}\). 
\item  \(X\) is a candidate rigid del Pezzo surface  with a basket of singularities \(\mathcal B\) if the multiplicities \(k_i\) are non-negative.
Repeat steps (iii)  to (vi) for each embedding computed at step (ii).
\end{enumerate}
\subsubsection{Pattern analysis}
In this step, we analyze the data of candidate examples obtained by using a computer algorithm to spot a pattern among the candidates to form a model of del Pezzo surfaces. The key indicator  of the pattern is usually the  number of orbifold points  and their multiplicities in the basket of orbifold points \(\mathcal B\).  
\subsubsection{Proving existence } The next step is to  prove the existence of such models as a complete intersection of some   $\wP$ or cone over it as a wellformed and quasismooth family of algebraic varieties.    The existence of a family of del Pezzo surfaces with  the  given Hilbert series follows straightforwardly from the existence of ambient \(\wP\). However, it does not guarantee that given del Pezzo surface has the same singularities as those predicted by the computer search algorithm of \ref{algorithm} or it is quasismooth. The singularities predicted by computer search may not necessarily lie on \(X\).  The crucial part of the existence is to show that  del Pezzo orbifold $X$ contains the same orbifold points as  suggested by the output in the computer search. 
  The wellformedness  can be deduced as  a consequence  of the existence of correct orbifold points, as  if \(X\) contains only point singularities then it avoids any  orbifold locus   of   $\P^6(w_{i})$ having a dimension greater than or equal to one.   

 The proof of the quasismoothness requires  an analysis of two different types of  loci:  orbifold loci and base loci.
 The orbifold loci are restrictions of the singular strata of \(\P^6(w_i) \) to \(X\)  and  the base  loci appear as  loci of the  linear systems of the successive intersection of weighted homogeneous forms while taking intersection with \(\wP\). A general member in each family del Pezzo surface appearing in these models is quasismooth, outside of these loci, due to    Bertini's theorem.
  
 \subsubsection{Non existence of degeneration to a toric variety}Each family in a given  models of orbifold del Pezzo surfaces is locally $\Q$-Gorenstein rigid. The first plurigenus   $h^0(-K_X)$ (and all  plurigenera $h^0(-mK))  $  are invariant under $\Q$-Gorenstein deformations. Moreover,   $h^0(-K_{X})>0$  for any toric Fano variety as  the  origin is always contained in the corresponding Fano polytope. Thus each family in a given model that satisfies    $h^0(-K_X)=0$  does not admit a $\Q$-Gorenstein degeneration to a toric variety. In total, 5 out of 13 models  have \(h^{0}(-K_X)=0\). These examples are of interest in the sense   they can not be constructed using usual toric methods.  

\subsection{Proof of Theorem \ref{Th:Index}}
We provide proof for one of the biregular  models appearing in Theorem \ref{Th:Index}. In all the  index models,  orbifold del Pezzo surfaces are obtained by using the intersection of some weighted \(\PI\) variety with a form of degree \(r\),  which can be assigned to the same variable in the
weight cube and in the equations  of \(\wP \) for each index model. Therefore we get the same base locus for each index model. Thus the below proof  presents base locus calculations for all models appearing in Theorem \ref{Th:Index}
and hence of quasismoothness.
    
\subsubsection{Biregular model \(\M\I_4\)} 
\label{S-Pf12}
This is an example of a  model where the base locus is bigger than zero dimensional but  geometrically it does not change for any family of the model, as described below. In this model, our parameter is \(r=3n+8\) for \(n\in \Z_{\ge 0}\), and the rest of the weights in terms of \(r\) are given as follows.   $$\begin{array}{lllll}p=r-3,& q=r-2,&s=r+2,&t=2r-3&m=2r-1.\end{array}$$  In this model each surface will have Fano index \(r\), so both the set of ambient weights and Fano index will vary for each  \(n\). If we choose a parameter \(w=(0,2:0,p:1,r)\), we get the embedding of weighted \(\PI\) 3-fold,   $$\wP \into \P(1,3,q,r^2,s,t,m) \text{ with  } K_{\wP}=\Oh(-2r).$$    Then  intersection of \(\wP\) with a generic weighted homogeneous form \(H_r\) of degree \(r\): 
$$X:= \wP\cap (H_r) \into \P(1,3,q,r,s,t,m)=\P(a,b,c,d,e,f,g)$$ is a  family of del Pezzo surfaces of Fano index  \(r\). The equations are presented by  $2 \times 2 $ minors of the weight cube: $$
 
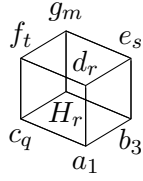
\begin{tikzpicture}[
  line join=round,
  y={(-0.86cm,0.36cm)},x={(.6cm,0.36cm)}, z={(0cm,.8cm)},
  arr/.style={line cap=round,shorten <= 1pt},baseline={([yshift=-.5em]current bounding box.center)}
]
\def\Side{1}
\coordinate (A1) at (0,0,0);\coordinate (A2) at (0,\Side,0);
\coordinate (A3) at (\Side,\Side,0);\coordinate (A4) at (\Side,0,0);\coordinate (B1) at (0,0,\Side);\coordinate (B2) at (0,\Side,\Side);\coordinate (B3) at (\Side,\Side,\Side);\coordinate (B4) at (\Side,0,\Side);

\draw[thin] (A2) -- (A1) -- (A4);\draw[thin] (B2) -- (B1) -- (B4) -- (B3) -- cycle;\draw[thin] (A1) -- (B1);\draw[thin] (A2) -- (B2);
\draw[thin] (A4) -- (B4);\draw[thin] (A2) -- (A3);
\draw[thin] (A3) -- (B3);\draw[thin] (A3) -- (A4);

\path[arr] 
  (A1) edge (A2)
  (B2) edge (A2)
  (B1) edge (B2)
  (B1) edge (A1)
  (B4) edge (A4)
  (B3) edge (A3)
  (B4) edge (B3)
  (A4) edge (A3);

\node[below] at (A1) {$a_1$};
\node[below] at (A2) {$c_q$};
\node[below] at (A3) {$H_{r}$};
\node[below] at (A4) {$b_3$};
\node[above] at (B1) {$d_r$};
\node[above] at (B2) {$f_t$};
\node[above] at (B3) {$g_m$};
\node[above] at (B4) {$e_s$};

\end{tikzpicture}$$
There are two different cases we need to discuss for the orbifold strata of \(X\):   \(r\) is odd and \(r\) is even.\\[2mm]
\textbf{Index \(r\) is odd}:  
 If \(X\) has an odd Fano index then all the weights of the ambient \(\P^6(w_i)\) are odd. None of the weights share any common factor except \(3, q\) and \(m\) as \(q\) and \(m\) are  multiples of \(3\). Thus it is clear that they are all at isolated orbifold points if they lie on  \(X\). The weight \(3\) locus is \[\wP\cap \P(3,q,m):=\V(bc,bg,cg)\subset \P(b,c,g). \] This locus is  manifestly 3 coordinate points \(b_{3}, c_q\) and \(g_m\). The coordinate points \(e_{s}\) and \(f_t\) also lies on \(X\). The coordinate point \(d_{r}\) does not lie on \(X\). 
 The summary of tangent and local variables on each coordinate point is given as follows.  
$$
\renewcommand*{\arraystretch}{1.4}
\begin{array}{cccll}
P_\orb & X\cap P_\orb  & \text {Tangent}|\text{local variables}  & \text{Conditions on forms}& \textrm{type }  \\\hline\hline
\evnrow        m & \text {coordinate pt } g_m & a,b,c,d\;|\;e,f & &\frac1m(s,t) \\
                t & \text {coordinate pt } f_t & a,b,e,d\;|\;c,m & H_r=d+\cdots &\frac1t(2,q)\\
\evnrow                                s & \text {coordinate pt-} e_s & a,c,d,f\;|\;b,g & H_r=d+\cdots&\frac1s(3,p) \\
                q & \text {coordinate pt } c_q & a,b,d,e\;|\;f,g & H_r=d+\cdots&\frac1q(1,1) \\
\evnrow                                3 & \text {coordinate pt-} b_3 & c,d,f,g\;|\;a,e & &\frac13(1,1)\\\hline

        \end{array}
$$
Therefore \(X\) is a family of del Pezzo surfaces with the desired orbifold points. Now we discuss the second case. \\[2mm] 
\textbf{Index \(r\) is even:} If the index \(r\) is even then the weights are \(q,r,s\) even and \(2\) is the common factor. Therefore we need to analyze the orbifold locus of weight \(2\) coming from the ambient \(\P^6(w_i)\). The rest of the orbifold strata is the same as in the case of an odd index. The weight two locus is  \[X\cap \P(q,r,s):=\V(d_rH_r,c_qe_s,H_re_s,d_rc_q)\subset \P(c,d,e). \] 
This is just a union of two coordinate points of  weight \(q\) and \(s\), which we already accounted for as isolated points. So we do not get any new orbifold points on this weight 2 locus.  

Now we discuss the base locus of the linear system \(|\Oh(r)|\).  None of the weights except \(1\) divides \(r\) so the linear system has base locus  \(\P(3,q,s,t,m)\). Then by using the equations of \(X\) the intersection of this locus  is the union of four copies of \(\P^1s\): 
\begin{equation}\label{baselocus}\P[e,g]\cup \P[b,e]\cup\P[f,g]\cup\P[c,f].\end{equation}Thus $X$ is quasismooth beside the reduced part of the above  locus and essentially this  locus does not change for any value of parameter \(r\). Therefore proving the quasismoothness for the first few values of \(r\) by using the computer algebra suffices the quasismoothness of \(\M\I_4\).  One can use the following  \textsc{Magma} code to establish the quasismoothness for any  $r(n)$.   
\begin{verbatim}
  procedure MI4(n)
    r:=3*n+8;
    q:=r-2;s:=r+2; t:=2*r-3;m:=2*r-1;
    P<x1,x2,x3,x4,x5,x6,x8>:=ProjectiveSpace(Rationals(),[1,r,q,t,3,s,m]);
    fr:=&+[ Random([1..5])*m : m in MonomialsOfWeightedDegree(CoordinateRing(P),r)];
    Eq:=[-x6*fr+x5*x8,-x4*fr+x3*x8,-x4*x6+x2*x8,-x4*x5+x1*x8,-x4*x5+x2*fr,
    -x4*x5+x3*x6,-x3*x5+x1*fr,-x2*x5+x1*x6,-x2*x3+x1*x4];
    X:=Scheme(P,Eq);
    SXred := ReducedSubscheme(JacobianSubrankScheme(X));
      if Dimension(SXred) eq -1 then 
      printf  "X is a quasismooth";
      end if; 
  end procedure;\end{verbatim} 

\subsection{Proof of Theorem \ref{Th:Weights}}
We replenish the proof for one  of the models appearing in Theorem \ref{Th:Weights}. The first step is to establish the  existence of models with specified orbifold points and numerical  invariants.   We assume that  a general member  $X$ of a family of del Pezzo surfaces has an embedding in  $ \P(a,b,c,d,e,f)$. We specify one of the weights to be \(r(n), n \in \Z^+\) and write the rest of the weights in terms of \(r\), given by   $$\begin{array}{lll} s=r+1,&u=r+3,&z=2r-1.\end{array}$$
\subsubsection{ Biregular Model $\M_{12}$} 
  This is an example of a model where the base locus is a finite number of points. We choose the  input parameter \(w=(0,1: 0,q:1,r)\) for \(r=2n\). This leads to  the  embedding of weighted \(\PI\) variety  $$\wP \into \P(2,s^{3},u^3,z) \text{ with  } K_{\wP}=\Oh(-3r-1).$$  We take an orbifold cone of weight \(r\) to get a Fano 4-fold \[\cC^4\wP\into\P(2,r,s^3,u^3,z)\] The cone variable of weight \(r\) will add \(-r\) to the canonical bundle. Then the intersection \[X=\cC^4\wP\cap (H_u)\cap (J_u)\into \P^6(2,r,s^3,u,z)=\P(a,b,c,d,e,f,g)\] with two forms of degrees \(u\) is a family of del Pezzo surfaces with \[K_X=\Oh(-(4r+1)+2u)=\Oh(-1).\] The extra cone variable \(r\) is not involved in the original equations of \(\wP\) but it appears in homogeneous form of degree \(u\)\(:\) \(H_u\) and \(J_u\). The equations of \(X\) are given by the  \(2\times 2 \) minors of the following cube.   
$$
 
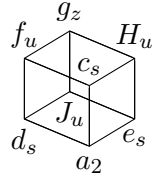
\begin{tikzpicture}[
  line join=round,
  y={(-0.86cm,0.36cm)},x={(.6cm,0.36cm)}, z={(0cm,.8cm)},
  arr/.style={line cap=round,shorten <= 1pt},baseline={([yshift=-.5em]current bounding box.center)}
]
\def\Side{1}
\coordinate (A1) at (0,0,0);\coordinate (A2) at (0,\Side,0);
\coordinate (A3) at (\Side,\Side,0);\coordinate (A4) at (\Side,0,0);\coordinate (B1) at (0,0,\Side);\coordinate (B2) at (0,\Side,\Side);\coordinate (B3) at (\Side,\Side,\Side);\coordinate (B4) at (\Side,0,\Side);

\draw[thin] (A2) -- (A1) -- (A4);\draw[thin] (B2) -- (B1) -- (B4) -- (B3) -- cycle;\draw[thin] (A1) -- (B1);\draw[thin] (A2) -- (B2);
\draw[thin] (A4) -- (B4);\draw[thin] (A2) -- (A3);
\draw[thin] (A3) -- (B3);\draw[thin] (A3) -- (A4);

\path[arr] 
  (A1) edge (A2)
  (B2) edge (A2)
  (B1) edge (B2)
  (B1) edge (A1)
  (B4) edge (A4)
  (B3) edge (A3)
  (B4) edge (B3)
  (A4) edge (A3);

\node[below] at (A1) {$a_{2}$};
\node[below] at (A2) {$d_s$};
\node[below] at (A3) {$J_{u}$};
\node[below] at (A4) {$e_s$};
\node[above] at (B1) {$c_s$};
\node[above] at (B2) {$f_u$};
\node[above] at (B3) {$g_z$};
\node[above] at (B4) {$H_u$};

\end{tikzpicture}$$

Now we check the singularities of  $ X$ on the orbifold locus of \(\P^6(w_i)\). The weight \(z\) locus is a coordinate point \(p_z\) which lies on \(X\).   The variables $a, c, d$  and \(e\) are tangent variables, and $b,f$ as local variables near  the coordinate  $g\ne 0$. Therefore it is an orbifold point of type   $\frac 1z(r,u)$. The weight \(u\) point does not lie on \(X\): as the weight \(u\) locus on \(\wP\) is the collection of three  coordinates points \(\P^2(u^3)\) and its intersection with two forms of weight \(u\) is an empty set.  The weight \(s\) locus   \[\P^2(c,d,e)= \P^2(s^3)\cap X\] is union of three coordinate points of \(\P^2(c,d,e)\).  One can verify using the implicit function theorem that  they are tree points of type \(\frac 1s(2,r)\). The weight \(r\) locus consists of \(\cC^4\wP\cap \P^3(r,u^3)\) as \(u=2r\) but it doesn't lie on \(X\). Therefore \(X\)  the basket \(\cB\) of orbifold points is \[\cB=\{3\times \frac1s(2,r), \frac1z(r,u)\}.\] 
 The base locus of the linear   system $|\Oh(u)|=<2^r,r^2,u>$ consists of \(\P(s^3,z)\) which intersects with \(X\) in exactly four points that consist
of  the basket \(\cB\) of \(X\) and they are  already proven to be quasismooth. Therefore \(X\) is quasismooth on the base locus as well. Thus  \(\M_{12} \) is a wellformed and quasismooth model of rigid orbifold del Pezzo surfaces. 

\section{Computer search  and sporadic examples}
\label{S-Sporadic}
\subsection{computational  output data }

\label{computer-search}
We use the algorithmic approach of \cite{QJSC, BKZ} to search for the candidate varieties for each Fano index \(I\). There is no finishing condition on the algorithm that provides \ the lists of all possible ambient weighted \(\PI\)  varieties, and thus its all possible quasilinear sections. However, our classification is complete   under   some  conditions,   chosen primarily due to the computational expensiveness of the search for the larger values of the ambient weights. We searched for all the cases  of  Fano index \(1\le I\le 16\)  and  adjunction number \(q\le 96\), i.e. we search for all possible del Pezzo surfaces \(X \subset \P^6(w_i)\) such that    \( \sum w_i+I \le 96   \).  \footnote{The detailed computer output can be accessed at https://sites.google.com/view/miqureshi/research/delPezzos. }.
\begin{rmk} We observe the following notable properties of our data.   \begin{enumerate}
\item 
 For the cases of index \(I\) greater than or equal to 5,  rigid   candidate surfaces  appear only for adjunction number \(6I.\)
\item We do not get any  rigid candidate examples of index 3 and 4. We do get candidate  examples with \(T\)-singularities.     
\item We only list those models for which at least two candidate examples appear in our computer search subject to the above given  bounds. Every example in  index \(\ge 2\) appears to be part of some model. The sporadic examples of the higher index are those for which we do not have two examples in their respective models, so we list them in the table of sporadic cases.

\item In index 1 and 16, there are sporadic examples, which are listed in Table \ref{Tab:Sporadic}. However, we expect the examples of index 16 to be part of some biregular model.  

\item In the index 2 cases, all working examples belong to some weight model. This is exactly similar to
the case of odPs in \(\PxP\) format of  \cite{QMOC}. Moreover 3 out of the 5 index 2, candidate models
do not work as quasismooth orbifolds and have the  worst singularities.

\item The candidates of an index greater than or equal to 3  seem to be  part of some index
model.
\item The model \(\M_{13}\) has the same numerics as a \( \PxP\) model  \(\boldsymbol P_{14}\) appeared  in \cite{QMOC}, i.e. They have the same ambient space, Hilbert series, degree and basket of orbifold points. \end{enumerate}
 \end{rmk}

We use the computer algebra system \textsc{Magma} to search for  the candidate examples and prove the quasismoothness of these models. The software  {\tt Mathematica} was used to calculate the generic formula for the \(-K_{X}^2\) for each model.
In Table \ref{tab-candidates},  \#Candidates denotes the total number of candidates obtained by using computer algebra and the row \#Examples lists the total number of wellformed and quasismooth working examples. The rest of the candidates obtained do not exist in \(\PI\)  format as they 
 contain
worst singularities.   \begin{table}[h]
\caption{Summary of computer search results for adjunction number \( \le 96.\)}\label{tab-candidates}\[
\renewcommand*{\arraystretch}{1.4}\begin{array}{c||cccccccccccccccc}
 \evnrow \textrm{Fano index}&1&2&3& 4 & 5 & 6 & 7 & 8 & 9 & 10 & 11 & 12 & 13 & 14 & 15 & 16    \\\hline
\#  \textrm{Candidates } & 33 & 18 & 0 & 0 & 1 & 0 & 1 & 1 & 2 & 1 & 3 & 1 & 3 & 5 & 2 & 5  \\
\evnrow \# \textrm{Examples }&25&9&0&0&1&0&1&1&2&1&3&1&3&5&2&5
\end{array}
\]

\end{table}
 
\subsection{Sporadic cases}
\label{Sporadic}
The following table lists all the sporadic examples. 
 \renewcommand*{\arraystretch}{1.2}
\begin{longtable}{>{\hspace{0.5em}}llccccr<{\hspace{0.5em}}}
\caption{Sporadic families of  rigid del  Pezzo surface of Fano index 1 and 16 in  $\PI $ format. } \label{Tab:Sporadic}\\
\toprule
S.No.&\multicolumn{1}{c}{WPS \& Para\ }&\multicolumn{1}{c}{Weight Cube}&\multicolumn{1}{c}{$\mathcal{B}$}&\multicolumn{1}{c}{$-K^2$}&\multicolumn{1}{c}{$h^0(-K)$}&\multicolumn{1}{c}{$I$}\\
\cmidrule(lr){1-1}\cmidrule(lr){2-2}\cmidrule(lr){3-3}\cmidrule(lr){4-4}\cmidrule(lr){5-7}
\endfirsthead
\multicolumn{7}{l}{\vspace{-0.25em}\scriptsize\emph{\tablename\ \thetable{} continued from previous page}}\\
\midrule
\endhead
\multicolumn{7}{r}{\scriptsize\emph{Continued on next page}}\\
\endfoot
\bottomrule
\endlastfoot
\evnrow  
 $\rownumber$& $\begin{matrix}
  \P(1^{7})\\w=(0,0:0,0: 1,1)
  \end{matrix} $&
 \begin{tikzpicture}[
  line join=round,
  y={(-0.86cm,0.36cm)},x={(.6cm,0.36cm)}, z={(0cm,.8cm)},
  arr/.style={line cap=round,shorten <= 1pt},baseline={([yshift=-.5em]current bounding box.center)}
]
\def\Side{1}
\coordinate (A1) at (0,0,0);\coordinate (A2) at (0,\Side,0);
\coordinate (A3) at (\Side,\Side,0);\coordinate (A4) at (\Side,0,0);\coordinate (B1) at (0,0,\Side);\coordinate (B2) at (0,\Side,\Side);\coordinate (B3) at (\Side,\Side,\Side);\coordinate (B4) at (\Side,0,\Side);

\draw[thin] (A2) -- (A1) -- (A4);\draw[thin] (B2) -- (B1) -- (B4) -- (B3) -- cycle;\draw[thin] (A1) -- (B1);\draw[thin] (A2) -- (B2);
\draw[thin] (A4) -- (B4);\draw[thin] (A2) -- (A3);
\draw[thin] (A3) -- (B3);\draw[thin] (A3) -- (A4);

\path[arr] 
  (A1) edge (A2)
  (B2) edge (A2)
  (B1) edge (B2)
  (B1) edge (A1)
  (B4) edge (A4)
  (B3) edge (A3)
  (B4) edge (B3)
  (A4) edge (A3);

\node[below] at (A1) {$1$};
\node[below] at (A2) {$1$};
\node[below] at (A3) {$1$};
\node[below] at (A4) {$1$};
\node[above] at (B1) {$1$};
\node[above] at (B2) {$1$};
\node[above] at (B3) {$\boldsymbol 1$};
\node[above] at (B4) {$1$};

\end{tikzpicture} &
  &$ 6         $&$ 7$ &1\\

\oddrow $\rownumber$&$ 
  \begin{matrix}
  \P(1,2^{3},3^{3})\\  
  w=(0,1:0,1: 1,2)
  \end{matrix}$&  
 \begin{tikzpicture}[
  line join=round,
  y={(-0.86cm,0.36cm)},x={(.6cm,0.36cm)}, z={(0cm,.8cm)},
  arr/.style={line cap=round,shorten <= 1pt},baseline={([yshift=-.5em]current bounding box.center)}
]
\def\Side{1}
\coordinate (A1) at (0,0,0);\coordinate (A2) at (0,\Side,0);
\coordinate (A3) at (\Side,\Side,0);\coordinate (A4) at (\Side,0,0);\coordinate (B1) at (0,0,\Side);\coordinate (B2) at (0,\Side,\Side);\coordinate (B3) at (\Side,\Side,\Side);\coordinate (B4) at (\Side,0,\Side);

\draw[thin] (A2) -- (A1) -- (A4);\draw[thin] (B2) -- (B1) -- (B4) -- (B3) -- cycle;\draw[thin] (A1) -- (B1);\draw[thin] (A2) -- (B2);
\draw[thin] (A4) -- (B4);\draw[thin] (A2) -- (A3);
\draw[thin] (A3) -- (B3);\draw[thin] (A3) -- (A4);

\path[arr] 
  (A1) edge (A2)
  (B2) edge (A2)
  (B1) edge (B2)
  (B1) edge (A1)
  (B4) edge (A4)
  (B3) edge (A3)
  (B4) edge (B3)
  (A4) edge (A3);

\node[below] at (A1) {$1$};
\node[below] at (A2) {$2$};
\node[below] at (A3) {$3$};
\node[below] at (A4) {$2$};
\node[above] at (B1) {$2$};
\node[above] at (B2) {$3$};
\node[above] at (B3) {$\boldsymbol 4$};
\node[above] at (B4) {$3$};

\end{tikzpicture} &$
3\times\frac{1}{3}(1,1)\cite{dp-CH}
$&1
 &$ 1$&1\\
 \evnrow  
 $\rownumber$& $\begin{matrix}
  \P(2,3^{3},4,5^2)\\w=(0,1:0,1: 2,4)
  \end{matrix} $&
 \begin{tikzpicture}[
  line join=round,
  y={(-0.86cm,0.36cm)},x={(.6cm,0.36cm)}, z={(0cm,.8cm)},
  arr/.style={line cap=round,shorten <= 1pt},baseline={([yshift=-.5em]current bounding box.center)}
]
\def\Side{1}
\coordinate (A1) at (0,0,0);\coordinate (A2) at (0,\Side,0);
\coordinate (A3) at (\Side,\Side,0);\coordinate (A4) at (\Side,0,0);\coordinate (B1) at (0,0,\Side);\coordinate (B2) at (0,\Side,\Side);\coordinate (B3) at (\Side,\Side,\Side);\coordinate (B4) at (\Side,0,\Side);

\draw[thin] (A2) -- (A1) -- (A4);\draw[thin] (B2) -- (B1) -- (B4) -- (B3) -- cycle;\draw[thin] (A1) -- (B1);\draw[thin] (A2) -- (B2);
\draw[thin] (A4) -- (B4);\draw[thin] (A2) -- (A3);
\draw[thin] (A3) -- (B3);\draw[thin] (A3) -- (A4);

\path[arr] 
  (A1) edge (A2)
  (B2) edge (A2)
  (B1) edge (B2)
  (B1) edge (A1)
  (B4) edge (A4)
  (B3) edge (A3)
  (B4) edge (B3)
  (A4) edge (A3);

\node[below] at (A1) {$2$};
\node[below] at (A2) {$4$};
\node[below] at (A3) {$5$};
\node[below] at (A4) {$3$};
\node[above] at (B1) {$3$};
\node[above] at (B2) {$5$};
\node[above] at (B3) {$\boldsymbol 6$};
\node[above] at (B4) {$\boldsymbol 4$};

\end{tikzpicture} &$5 \times \frac13(1,1),2\times \frac15(1,1)$
  &$ \frac 4{15}       $&$ 0$ &$ 1$\\

\oddrow $\rownumber$&$ 
  \begin{matrix}
  \P(1,3^{3},5^{3})\\  
  w=(0,2:0,2: 1,3)
  \end{matrix}$&  
 \begin{tikzpicture}[
  line join=round,
  y={(-0.86cm,0.36cm)},x={(.6cm,0.36cm)}, z={(0cm,.8cm)},
  arr/.style={line cap=round,shorten <= 1pt},baseline={([yshift=-.5em]current bounding box.center)}
]
\def\Side{1}
\coordinate (A1) at (0,0,0);\coordinate (A2) at (0,\Side,0);
\coordinate (A3) at (\Side,\Side,0);\coordinate (A4) at (\Side,0,0);\coordinate (B1) at (0,0,\Side);\coordinate (B2) at (0,\Side,\Side);\coordinate (B3) at (\Side,\Side,\Side);\coordinate (B4) at (\Side,0,\Side);

\draw[thin] (A2) -- (A1) -- (A4);\draw[thin] (B2) -- (B1) -- (B4) -- (B3) -- cycle;\draw[thin] (A1) -- (B1);\draw[thin] (A2) -- (B2);
\draw[thin] (A4) -- (B4);\draw[thin] (A2) -- (A3);
\draw[thin] (A3) -- (B3);\draw[thin] (A3) -- (A4);

\path[arr] 
  (A1) edge (A2)
  (B2) edge (A2)
  (B1) edge (B2)
  (B1) edge (A1)
  (B4) edge (A4)
  (B3) edge (A3)
  (B4) edge (B3)
  (A4) edge (A3);

\node[below] at (A1) {$1$};
\node[below] at (A2) {$3$};
\node[below] at (A3) {$5$};
\node[below] at (A4) {$3$};
\node[above] at (B1) {$3$};
\node[above] at (B2) {$5$};
\node[above] at (B3) {$\boldsymbol 7$};
\node[above] at (B4) {$5$};

\end{tikzpicture} &$
3\times\frac{1}{3}(1,1), 3\times \frac15(1,1)
$&$\frac25
 $&$ 1$&1\\
 
 \evnrow  
 $\rownumber$& $\begin{matrix}
  \P(2,3,4^2,5^3)\\w=(0,1:0,1: 3,4)
  \end{matrix}$&
 \begin{tikzpicture}[
  line join=round,
  y={(-0.86cm,0.36cm)},x={(.6cm,0.36cm)}, z={(0cm,.8cm)},
  arr/.style={line cap=round,shorten <= 1pt},baseline={([yshift=-.5em]current bounding box.center)}
]
\def\Side{1}
\coordinate (A1) at (0,0,0);\coordinate (A2) at (0,\Side,0);
\coordinate (A3) at (\Side,\Side,0);\coordinate (A4) at (\Side,0,0);\coordinate (B1) at (0,0,\Side);\coordinate (B2) at (0,\Side,\Side);\coordinate (B3) at (\Side,\Side,\Side);\coordinate (B4) at (\Side,0,\Side);

\draw[thin] (A2) -- (A1) -- (A4);\draw[thin] (B2) -- (B1) -- (B4) -- (B3) -- cycle;\draw[thin] (A1) -- (B1);\draw[thin] (A2) -- (B2);
\draw[thin] (A4) -- (B4);\draw[thin] (A2) -- (A3);
\draw[thin] (A3) -- (B3);\draw[thin] (A3) -- (A4);

\path[arr] 
  (A1) edge (A2)
  (B2) edge (A2)
  (B1) edge (B2)
  (B1) edge (A1)
  (B4) edge (A4)
  (B3) edge (A3)
  (B4) edge (B3)
  (A4) edge (A3);

\node[below] at (A1) {$3$};
\node[below] at (A2) {$4$};
\node[below] at (A3) {$5$};
\node[below] at (A4) {$\boldsymbol 4$};
\node[above] at (B1) {$4$};
\node[above] at (B2) {$5$};
\node[above] at (B3) {$\boldsymbol 6$};
\node[above] at (B4) {$5$};

\end{tikzpicture} &$3\times \frac15(1,2)$
  &$ \frac {1}{5}       $&$ 0$&1 \\

\oddrow $\rownumber$&$ 
  \begin{matrix}
  \P(3^{2},5^{2},7,9,11)\\  
  w=(0,2:0,2: 2,7)
  \end{matrix}$&  
 \begin{tikzpicture}[
  line join=round,
  y={(-0.86cm,0.36cm)},x={(.6cm,0.36cm)}, z={(0cm,.8cm)},
  arr/.style={line cap=round,shorten <= 1pt},baseline={([yshift=-.5em]current bounding box.center)}
]
\def\Side{1}
\coordinate (A1) at (0,0,0);\coordinate (A2) at (0,\Side,0);
\coordinate (A3) at (\Side,\Side,0);\coordinate (A4) at (\Side,0,0);\coordinate (B1) at (0,0,\Side);\coordinate (B2) at (0,\Side,\Side);\coordinate (B3) at (\Side,\Side,\Side);\coordinate (B4) at (\Side,0,\Side);

\draw[thin] (A2) -- (A1) -- (A4);\draw[thin] (B2) -- (B1) -- (B4) -- (B3) -- cycle;\draw[thin] (A1) -- (B1);\draw[thin] (A2) -- (B2);
\draw[thin] (A4) -- (B4);\draw[thin] (A2) -- (A3);
\draw[thin] (A3) -- (B3);\draw[thin] (A3) -- (A4);

\path[arr] 
  (A1) edge (A2)
  (B2) edge (A2)
  (B1) edge (B2)
  (B1) edge (A1)
  (B4) edge (A4)
  (B3) edge (A3)
  (B4) edge (B3)
  (A4) edge (A3);

\node[below] at (A1) {$3$};
\node[below] at (A2) {$7$};
\node[below] at (A3) {$\boldsymbol 9$};
\node[below] at (A4) {$5$};
\node[above] at (B1) {$5$};
\node[above] at (B2) {$9$};
\node[above] at (B3) {$11$};
\node[above] at (B4) {$\boldsymbol 7$};

\end{tikzpicture} &$
\begin{matrix}5\times\frac{1}{3}(1,1), 2\times \frac15(1,1),\\\frac1{11}(1,3)\end{matrix}
$&$\frac{14}{165}
 $&$ 0$&1\\
 \evnrow  
 $\rownumber$& $\begin{matrix}
  \P(1,5^2,7,9,11^2)\\w=(0,4:0,4: 1,7)
  \end{matrix}$&
 \begin{tikzpicture}[
  line join=round,
  y={(-0.86cm,0.36cm)},x={(.6cm,0.36cm)}, z={(0cm,.8cm)},
  arr/.style={line cap=round,shorten <= 1pt},baseline={([yshift=-.5em]current bounding box.center)}
]
\def\Side{1}
\coordinate (A1) at (0,0,0);\coordinate (A2) at (0,\Side,0);
\coordinate (A3) at (\Side,\Side,0);\coordinate (A4) at (\Side,0,0);\coordinate (B1) at (0,0,\Side);\coordinate (B2) at (0,\Side,\Side);\coordinate (B3) at (\Side,\Side,\Side);\coordinate (B4) at (\Side,0,\Side);

\draw[thin] (A2) -- (A1) -- (A4);\draw[thin] (B2) -- (B1) -- (B4) -- (B3) -- cycle;\draw[thin] (A1) -- (B1);\draw[thin] (A2) -- (B2);
\draw[thin] (A4) -- (B4);\draw[thin] (A2) -- (A3);
\draw[thin] (A3) -- (B3);\draw[thin] (A3) -- (A4);

\path[arr] 
  (A1) edge (A2)
  (B2) edge (A2)
  (B1) edge (B2)
  (B1) edge (A1)
  (B4) edge (A4)
  (B3) edge (A3)
  (B4) edge (B3)
  (A4) edge (A3);

\node[below] at (A1) {$1$};
\node[below] at (A2) {$7$};
\node[below] at (A3) {$11$};
\node[below] at (A4) {$5$};
\node[above] at (B1) {$5$};
\node[above] at (B2) {$11$};
\node[above] at (B3) {$\boldsymbol{15}$};
\node[above] at (B4) {$9$};

\end{tikzpicture} &$\begin{matrix}\frac{1}{7}(1,1),  \frac19(1,1),\\2\times \frac1{11}(1,7)\end{matrix}$
  &$ \frac {74}{693}       $&$ 1$ &1\\

\oddrow $\rownumber$&$ 
  \begin{matrix}
  \P(1,6,8,12,13,17,19)\\  
  w=(0,5:0,7: 1,12)
  \end{matrix}$& 
 \begin{tikzpicture}[
  line join=round,
  y={(-0.86cm,0.36cm)},x={(.6cm,0.36cm)}, z={(0cm,.8cm)},
  arr/.style={line cap=round,shorten <= 1pt},baseline={([yshift=-.5em]current bounding box.center)}
]
\def\Side{1}
\coordinate (A1) at (0,0,0);\coordinate (A2) at (0,\Side,0);
\coordinate (A3) at (\Side,\Side,0);\coordinate (A4) at (\Side,0,0);\coordinate (B1) at (0,0,\Side);\coordinate (B2) at (0,\Side,\Side);\coordinate (B3) at (\Side,\Side,\Side);\coordinate (B4) at (\Side,0,\Side);

\draw[thin] (A2) -- (A1) -- (A4);\draw[thin] (B2) -- (B1) -- (B4) -- (B3) -- cycle;\draw[thin] (A1) -- (B1);\draw[thin] (A2) -- (B2);
\draw[thin] (A4) -- (B4);\draw[thin] (A2) -- (A3);
\draw[thin] (A3) -- (B3);\draw[thin] (A3) -- (A4);

\path[arr] 
  (A1) edge (A2)
  (B2) edge (A2)
  (B1) edge (B2)
  (B1) edge (A1)
  (B4) edge (A4)
  (B3) edge (A3)
  (B4) edge (B3)
  (A4) edge (A3);

\node[below] at (A1) {$1$};
\node[below] at (A2) {$12$};
\node[below] at (A3) {$17$};
\node[below] at (A4) {$6$};
\node[above] at (B1) {$8$};
\node[above] at (B2) {$19$};
\node[above] at (B3) {$\boldsymbol{24}$};
\node[above] at (B4) {$13$};

\end{tikzpicture} &$
\begin{matrix}\frac{1}{13}(1,4),  \frac1{17}(1,2),\\\frac1{19}(1,7)\end{matrix}
$&$\frac{202}{4199}
 $&$ 1$&1\\\hline\hline
 \evnrow  
 $\rownumber$& $\begin{matrix}
   \P(3,5,14,16,27,28,29)\\  
  w=(0,2:0,11: 3,16)
  \end{matrix}$&
 \begin{tikzpicture}[
  line join=round,
  y={(-0.86cm,0.36cm)},x={(.6cm,0.36cm)}, z={(0cm,.8cm)},
  arr/.style={line cap=round,shorten <= 1pt},baseline={([yshift=-.5em]current bounding box.center)}
]
\def\Side{1}
\coordinate (A1) at (0,0,0);\coordinate (A2) at (0,\Side,0);
\coordinate (A3) at (\Side,\Side,0);\coordinate (A4) at (\Side,0,0);\coordinate (B1) at (0,0,\Side);\coordinate (B2) at (0,\Side,\Side);\coordinate (B3) at (\Side,\Side,\Side);\coordinate (B4) at (\Side,0,\Side);

\draw[thin] (A2) -- (A1) -- (A4);\draw[thin] (B2) -- (B1) -- (B4) -- (B3) -- cycle;\draw[thin] (A1) -- (B1);\draw[thin] (A2) -- (B2);
\draw[thin] (A4) -- (B4);\draw[thin] (A2) -- (A3);
\draw[thin] (A3) -- (B3);\draw[thin] (A3) -- (A4);

\path[arr] 
  (A1) edge (A2)
  (B2) edge (A2)
  (B1) edge (B2)
  (B1) edge (A1)
  (B4) edge (A4)
  (B3) edge (A3)
  (B4) edge (B3)
  (A4) edge (A3);

\node[below] at (A1) {$3$};
\node[below] at (A2) {$16$};
\node[below] at (A3) {$18$};
\node[below] at (A4) {$5$};
\node[above] at (B1) {$14$};
\node[above] at (B2) {$27$};
\node[above] at (B3) {$29$};
\node[above] at (B4) {$\boldsymbol{16}$};

\end{tikzpicture} &$\begin{matrix}\frac{1}{3}(1,1),  \frac15(1,1), \frac1{14}(1,9),\\\frac1{27}(1,4), \frac1{29}(1,16) \end{matrix}$
  &$ \frac {47104}{27405}       $&$ 2$ &16\\

\oddrow $\rownumber$&$ 
  \begin{matrix}
  \P(7,10,13,16,19,22,25)\\  
  w=(0,3:0,6: 7,16)
  \end{matrix}$& 
 \begin{tikzpicture}[
  line join=round,
  y={(-0.86cm,0.36cm)},x={(.6cm,0.36cm)}, z={(0cm,.8cm)},
  arr/.style={line cap=round,shorten <= 1pt},baseline={([yshift=-.5em]current bounding box.center)}
]
\def\Side{1}
\coordinate (A1) at (0,0,0);\coordinate (A2) at (0,\Side,0);
\coordinate (A3) at (\Side,\Side,0);\coordinate (A4) at (\Side,0,0);\coordinate (B1) at (0,0,\Side);\coordinate (B2) at (0,\Side,\Side);\coordinate (B3) at (\Side,\Side,\Side);\coordinate (B4) at (\Side,0,\Side);

\draw[thin] (A2) -- (A1) -- (A4);\draw[thin] (B2) -- (B1) -- (B4) -- (B3) -- cycle;\draw[thin] (A1) -- (B1);\draw[thin] (A2) -- (B2);
\draw[thin] (A4) -- (B4);\draw[thin] (A2) -- (A3);
\draw[thin] (A3) -- (B3);\draw[thin] (A3) -- (A4);

\path[arr] 
  (A1) edge (A2)
  (B2) edge (A2)
  (B1) edge (B2)
  (B1) edge (A1)
  (B4) edge (A4)
  (B3) edge (A3)
  (B4) edge (B3)
  (A4) edge (A3);

\node[below] at (A1) {$7$};
\node[below] at (A2) {$16$};
\node[below] at (A3) {$19$};
\node[below] at (A4) {$10$};
\node[above] at (B1) {$13$};
\node[above] at (B2) {$22$};
\node[above] at (B3) {$25$};
\node[above] at (B4) {$\boldsymbol{16}$};

\end{tikzpicture} &$
\begin{matrix}\frac{1}{7}(1,2),  \frac1{10}(1,3), \frac1{13}(1,5),\\\frac1{19}(1,8), \frac1{22}(1,7),\frac1{25}(1,2)\end{matrix}
$&$\frac{278768}{475475}
 $&$ 1$&16

\\
 \evnrow  
 $\rownumber$& $\begin{matrix}
  \P(6,11^2,16,21^2,26)\\w=(0,5:0,5: 6,16)
  \end{matrix}$&
 \begin{tikzpicture}[
  line join=round,
  y={(-0.86cm,0.36cm)},x={(.6cm,0.36cm)}, z={(0cm,.8cm)},
  arr/.style={line cap=round,shorten <= 1pt},baseline={([yshift=-.5em]current bounding box.center)}
]
\def\Side{1}
\coordinate (A1) at (0,0,0);\coordinate (A2) at (0,\Side,0);
\coordinate (A3) at (\Side,\Side,0);\coordinate (A4) at (\Side,0,0);\coordinate (B1) at (0,0,\Side);\coordinate (B2) at (0,\Side,\Side);\coordinate (B3) at (\Side,\Side,\Side);\coordinate (B4) at (\Side,0,\Side);

\draw[thin] (A2) -- (A1) -- (A4);\draw[thin] (B2) -- (B1) -- (B4) -- (B3) -- cycle;\draw[thin] (A1) -- (B1);\draw[thin] (A2) -- (B2);
\draw[thin] (A4) -- (B4);\draw[thin] (A2) -- (A3);
\draw[thin] (A3) -- (B3);\draw[thin] (A3) -- (A4);

\path[arr] 
  (A1) edge (A2)
  (B2) edge (A2)
  (B1) edge (B2)
  (B1) edge (A1)
  (B4) edge (A4)
  (B3) edge (A3)
  (B4) edge (B3)
  (A4) edge (A3);

\node[below] at (A1) {$6$};
\node[below] at (A2) {$16$};
\node[below] at (A3) {$21$};
\node[below] at (A4) {$11$};
\node[above] at (B1) {$11$};
\node[above] at (B2) {$21$};
\node[above] at (B3) {$26$};
\node[above] at (B4) {$\boldsymbol{16}$};

\end{tikzpicture} &$\begin{matrix}\frac{1}{6}(1,1),  2\times \frac1{11}(1,5),\\2\times \frac1{21}(1,10),\frac{1}{26}(1,1) \end{matrix}$
  &$ \frac {2048}{3003}       $&$ 1$ &16\\

\end{longtable}

%%%%%%%%%%%%%%%%%%%%%%%%%%%%%%%%%%%%

\subsection*{Acknowledgments}

I am thankful to  Alexander Kasprzyk for useful discussions and to
referees for their comments that helped improved the exposition of the article. I am also thankful to the anonymous referees for the careful reading of the
manuscript that lead to the improvement of the exposition of this article. I would  like to acknowledge the support provided by the IRC Intelligent and Secure System at the King Fahd University of
Petroleum and Minerals  through Project No. INSS2308. 

%%%%%%%%%%%%%%%%%%%%%%%%%%%%%%%%%%%%
\bibliographystyle{amsalpha}
\bibliography{References}
%%%%%%%%%%%%%%%%%%%%%%%%%%%%%%%%%%%%

\noindent {\sc Muhammad Imran Qureshi, \\
Department of Mathematics, King Fahd University of Petroleum and Minerals, Dhahran, 31261, Saudi Arabi \&\\ 
   Interdisciplinary Research Center for Intelligent and Secure Systems, King Fahd University of Petroleum and Minerals, Dhahran, 31261, Saudi Arabia }\\
\noindent {\tt Email: imran.qureshi@kfupm.edu.sa}
%%%%%%%%%%%%%%%%%%%%%%%%%%%%%%%%%%%%
\end{document}